 \documentclass[final,3p,times]{elsarticle}
\usepackage{amssymb}
\usepackage{graphicx,amsmath}
\usepackage{ifthen}
\usepackage{endnotes}
\usepackage{hyperref}
\hypersetup{colorlinks=true,citecolor=blue}
\def\beq{\begin{equation}}
\def\eeq{\end{equation}}
\def\bea{\begin{eqnarray}}
\def\eea{\end{eqnarray}}
\def\nn{\nonumber}

\biboptions{square,comma,numbers,sort&compress}

\journal{Technology Forecasting and Social Change}

\begin{document}

\begin{frontmatter}


\title{On the changeover timescales of technology transitions and induced efficiency changes: an overarching theory}%\tnoteref{label1}}
\author{Jean-Fran\c{c}ois Mercure \corref{cor1}}
\ead{jm801@cam.ac.uk}
%\ead[url]{www.e3mgmodel.com}
%%\fntext[label2]{}
\cortext[cor1]{Corresponding author: Jean-Fran\c{c}ois Mercure}
\address{Cambridge Centre for Climate Change Mitigation Research (4CMR), Department of Land Economy, University of Cambridge, 19 Silver Street, Cambridge, CB3 1EP, United Kingdom}

\begin{abstract}

This paper presents a general theory that aims at explaining timescales observed empirically in technology transitions and predicting those of future transitions. This framework is used further to derive a theory for exploring the dynamics that underlie the complex phenomenon of irreversible and path dependent price or policy induced efficiency changes. Technology transitions are known to follow patterns well described by logistic functions, which should more rigorously be modelled mathematically using the Lotka-Volterra family of differential equations, originally developed to described the population growth of competing species. The dynamic evolution of technology has also been described theoretically using evolutionary dynamics similar to that observed in nature. The theory presented here joins both approaches and presents a methodology for predicting changeover time constants in order to describe real systems of competing technologies. The problem of price or policy induced efficiency changes becomes naturally explained using this framework. Examples of application are given.

\end{abstract}

\begin{keyword}
Technological Change \sep Technology Transitions \sep Induced technological change

\end{keyword}

\end{frontmatter}

% \linenumbers

\section{Introduction}

\subsection{The challenge of describing induced technological change}

The challenge of predicting the future energy transition is not a new one, but one of considerable importance \cite{Wilson2011,Grubler2012,Grubler1999,Grubler1998,Azar1999,IPCCSRES}. This is due to the fact that there is significant uncertainty in the technical feasibility of a global technology transformation that would enable to reduce global greenhouse emissions to a small fraction of their current levels, in order to limit anthropogenic climate change while maintaining the accessibility of energy-related services and a healthy economic system. Changes in the prices of energy carriers have the potential to generate significant inflation and economic slowdown (e.g. the price of electricity or oil, for a review of the matter see Jones $et$ $al.$\cite{Jones2004}, and for a recent model by the author of energy commodity prices, see \cite{MercureSalas2012,MercureSalas2012b}), which could occur in some scenarios of energy technology transitions where, for instance, the price of carbon is passed on by utility companies into the price of electricity. Technology however does change endogenously, an aspect that has not been explored in a large amount of detail but is the subject of an emerging literature (see for instance \cite{Kohler2006, Grubb2006, Crassous2006, Edenhofer2006, Azar1999, Grubler1999,Grubler1998} and the Innovation Modelling Comparison Project \cite{IMCP2006}.). Within this approach, forecasting technology is increasingly seen as a key to sensible climate policy-making, in particular for exploring potential technology support policy portfolios and their impact on emissions.

Many current models of technology (bottom-up models) use exogenous (i.e. imposed) timescales and rates of technological change, which therefore do not depend on actual economic indicators. Meanwhile, many economic (top-down) models rely on a so-called Autonomous Energy Efficiency Indicator (AEEI), which extrapolates the global non price-related energy efficiency of the economy from historical data \cite{Azar1999,Kaufmann2004}, associated directly with technology improvements. As Dowlatabadi and Oravetz argue, "We need only delve into the recent past to realize that a formulation of technological change that ignores price effects may be unrealistic"\cite{Dowlatabadi2006}. Most formulations of the AEEI lead to reasoning that is nearly  circular, where one is trying to model large future changes in energy efficiency produced by policy using a historical trend where no such strong climate policy precedent exists, choosing as a result the assumption. While most models of global energy systems and greenhouse gas emissions possess highly detailed and accurate accounting of global emissions sources, significant uncertainty or inaccuracies arise in their forecasted evolution of emissions factors up to standard horizons of 2050 or 2100, given chosen policy portfolios, due to a lack of dynamic interactions between components in their representation of technology substitution (e.g. changeover timescales, technology learning, cross-sectoral synergies and interference, etc, for reviews of models see for instance Barker $et$ $al.$\cite{IPCCAR4Ch11} and Edenhofer $et$ $al.$ \cite{Edenhofer2010}), or even the simple absence of a backstop technology (for instance in \cite{Edenhofer2006}). In particular, no consensus exists on the value of the carbon price necessary to decarbonise the global economy, for which values between 100 and above 1500~\$/tCO$_2$ have been claimed to be necessary (see for instance \cite{Edenhofer2006}). Since the decarbonisation process through the use of a carbon price signal relies primarily on technological change and diffusion, the value of the carbon price depends heavily on the particular model chosen and the amount of substitution possibilities and dynamic interactions present in the model for representing technological change. This motivates strongly the development and expansion of theoretical and computational models of technology diffusion, for application in energy supply (e.g. power) as well as end-use sectors (e.g. appliances, transport, industry).

Insight on the structure of a future energy transition has been sought in historical data \citep{Grubler2012,Grubler1998}. A well established literature of empirical modelling of technological change exists, which provides a wealth of details on past evolutions of, for instance, energy systems, energy end-use appliances, innovation and the diffusion of new consumption goods. In particular, changes in the consumption of global primary energy has famously been shown to follow logistic trends \citep{Marchetti1978}. Early work by Bass presents an attempt to model the diffusion of home appliances using a logistic model \cite{Bass1969}. Fisher and Pry analysed the diffusion of 17 different pairwise substitutions of natural versus synthetic materials or other types in several sectors, all of which follow a simple logistic model \cite{Fisher1971}. The diffusion of infrastructure was shown to follow similar trends \citep{Grubler1999}, as well as the transition from horses to petrol cars in the 20's \citep{Nakicenovic1986}. Sharif and Kabir showed that the replacement of steam for electric and diesel locomotives in U.S. railroads, as well as metal for wood construction materials in the U.S. marine follow a special form of the logistic substitution model \cite{Sharif1976}. More recently, an in-depth analysis of power technologies was carried out by Wilson \cite{Wilson2009}. Many more such analyses have successfully applied one or another form of the logistic model to product or technology substitution in the marketplace \cite{Farrell1993,Grubler1998,Grubler1999}. These works provide important insight on the structure of endogenous changes occurring within all sectors of society that rely on technology and its improvement in time.

The suggestion that multi-technology\footnote{i.e. more than two} competition in the marketplace could follow patterns obeyed by competing species in ecosystems was explored by Barghava \cite{Bhargava1989}, who mapped the problem to the Lotka-Volterra family of coupled population growth equations for biological systems \citep{Lotka1925,Volterra1939}, themselves a generalisation of the classic logistic model of Verhulst \cite{Verhulst1838} for forecasting the French population after 1838.  Possible transitions in the transport sector were explored by Gr\"ubler \cite{Grubler1990}, where it was also suggested that the competition between technologies could follow the Lotka-Volterra family of equations. However, the nature of the exchange between technologies, as well as the nature of the timescales involved, or the parameterisation of the equations, has never been satisfactorily clarified in the literature. This concerns parameters that are related, but do not correspond exactly to, to birth rates and death rates in biological systems, but that give rise to transition time constants. Timescales involved in technological transitions have been measured empirically extensively, however no theoretical framework exists to the author's knowledge that can predict what they should be, even though a simple modeller's intuition tells us that it should involve birth and death rates --- technology construction and decommission timescales, and/or bottlenecks in production supply chain and their timescales of expansion. Additionally, confusion and disagreement exists between experts in the field as to whether the Lotka-Volterra family of equations represents or not simultaneous pairwise interactions between species (technologies) where exchanges occur between all competing types with one-another \cite{Grubler2012Comm}.

Following previous work where I introduced a model in which the Lotka-Volterra family of equations is used to model power sector technology substitution \cite{Mercure2012}, in this paper I attempt to provide a theoretical clarification of the process of technology competition and substitution, and to bring this framework to a useful generalisation that can be applied to various sectors of the global economy. In particular, I explore the origin and properties of the timescales involved in such technological change. A complete model for the derivation of changeover timescales is given, which belong to the interacting pair rather than to individual technologies. Since this mathematical problem is one of many-body interactions, it will be stressed that no analytical solution exists and that solutions should simply be calculated dynamically using discrete time-steps.

This process of population movements between elements of a given list of existing technologies with different properties underlies changes in technology properties averaged over the whole sector, such as the efficiency of energy use in transport, steel production, electricity production, lighting, etc. It is known to be asymmetric and path dependent \cite{Cologni2009}, and therefore hysteretic and multi-valued.\footnote{The response of the oil intensity of the global economy with respect to the oil shocks of the 70's-80's has been shown asymmetric \cite{Cologni2009}.} This means that while an increase in the price of energy generates an increase in the efficiency of energy use, a subsequent decrease in price does not produce a corresponding decrease in efficiency and the system does not return to its starting point, due to permanent changes of technology. This makes the use of linear correlations (e.g. econometric) as relationships for forecasting inappropriate or at least incomplete,\footnote{i.e. the parameters of such correlations are time dependent; they change when the historical period chosen is changed, due to permanent internal structure changes, e.g. permanent changes in technology.} because efficiency is thus not a single valued function of prices, gross domestic product and other economic variables or simply time. While the replacement of parameters such as the AEEI in economic models is a difficult task, it must involve (1) technical properties of individual technologies and (2) an appropriate dynamic weighing of these properties as the population of technologies evolves in time. 

Finally, the process of technological change has often been qualitatively compared with the theory of evolution. Theoretical work by Dercole $et$ $al.$ shows that technological change can follow evolutionary dynamics by creating a parallel in which technology characteristics endogenously improve in order to increase their market share \cite{Dercole2008}. Such dynamics increase the technology diversity that exists at any time. Changes in the economic/market/policy environment can change which technology (specie) strives best and give rise to new technologies diffusing massively and displacing other older systems. The key to such a process being the underlying diversity, such environmental changes can uncover existing technology potential which, before the changes, survived within small market niches provided by specialist applications. Such a framework is also consistent with observed real technology dynamics \cite{Wilson2011}.

The first half of this work follows a simple progression where the process of technological change is initially derived, for which six definitions and three assumptions are given, and terms are derived for, in order, exchanges of capacity between technologies, increases in capacity requirements and permanent decommissions of capacity. These terms are combined in order to construct the complete equation for technological change. This enables a detailed analysis of the parameterisation of logistic transition models and a determination of technology changeover time constants. This theoretical framework is then used to describe the time evolution of sectoral parameters such as energy efficiency or intensity. In particular, path dependent and hysteresis properties are analysed. Their implications on the understanding of the process are stressed. Finally, the parallel with evolutionary dynamics is discussed.

\section{The theoretical model}
\subsection{Competing markets and production}
Each individual technology within a competing set strive for gaining ground in market space, which is defined in terms of fractions of the total demand $D(t)$ for a particular common good or service. This includes for example fuels (electricity, ethanol, petrol, etc), materials such as steel or services such as transport, lighting, communication/connectivity, etc. For each of these, a set of substitutable technologies $i = 1, 2, ..., N$ compete and produce a service output $G_i(t)$. The units of $G_i$ are the key for its definition, and may correspond to, for example, energy (MWh, GJ, Mtoe...), mass (kg of materials), movement (person-km or ton-km), light (J/m$^2$), etc. 

The production of these services is performed by a production capacity $U_i$, which is used as the service is needed. This therefore invokes the definition of the intensity of capacity use, which we term the capacity factor $CF_i$, corresponding to the fraction of time for which the unit is producing at it's full capacity, or alternatively, the fraction of its total capacity for which it is being used, or a mixture of both.\footnote{No distinction is made between for instance having two units where only one is used continuously, several units all of which are used 50\% of the time, or several units used at 50\% of their capacity, or any other such combination giving $CF = 50\%$.} Given these definitions, the generation of a service $G_i$ is 
\beq
G_i(t) = CF_i(t)U_i(t),
\eeq
which sums up to the total demand $D = \sum_i G_i$, with a total capacity $U_{tot}$ and overall average capacity factor $\overline{CF}$:
\beq
D(t) = \overline{CF}(t) U_{tot}(t).
\eeq
The composition of the production $G_i$ is therefore dependent, following individual capacity factors $CF_i$, on the number of producing units in each technology category $U_i$, which defines the composition of the sector. This can be expressed as a fraction of the total, or shares of production capacity $S_i(t) = U_i/U_{tot}$. The capacity can be expressed in terms of the shares $S$ and the system weighed average capacity cactor $\overline{CF}$:
\bea
\label{eq:U}
U_i(t) &=& {S_i(t) D(t) \over \overline{CF}(t)},\\
\overline{CF}(t) &=& \sum_i S_i(t)CF_i(t),\nn
\eea
and the market shares (shares of demand) are
\beq
\sigma_i(t) = {CF_i(t) \over \overline{CF}(t)} S_i(t).
\eeq
These equations were already derived for the power sector in previous work \citep{Mercure2012}, in which a differential form is invoked for explanatory purpose,
\beq
dU_i =  {D \over \overline{CF}} dS_i + {S_i \over \overline{CF}} dD - {S_i D \over \overline{CF}^2}d\overline{CF},
\eeq
This equation demonstrates that the number of units of each type of technology can vary following three types of changes: changes in the technology mix $dS_i$ (the exchange term, where technologies exchange capacity shares), changes in the total demand $dD$ and changes in the overall efficiency over which the system is used $d\overline{CF}$. Additionally, however, the total number of units can be fixed ($dU_i=0, dS_i = 0$) and changes in demand can be compensated by changes in overall capacity factor within limits, the composition of the system could change at constant demand and efficiency ($dD = 0, d\overline{CF}=0$), or the number of units could change at constant demand and composition, compensated by changes in efficiency ($dS_i = 0, dD = 0$), etc. In general however, changes in $U_i$ originate from simultaneous changes in all of these independent variables. In particular, the exchange term is written terms of the sum of individual interactions between pairs of technologies:
\beq
dS_i = \sum_{j=1}^N dS_{ij} = \sum_{j=1}^N {1 \over T_{ij}} S_i S_j dt,
\label{eq:LotkaVolterra}
\eeq
where changes in shares during a time interval $dt$ occur due to (1) the growth rate of a technology, which is proportional to its own extent of diffusion (or share of the total capacity) $S_i$, as well as (2) to the rate at which it can replace other technologies, events that occur every time units are decommissioned, where the decommission rates are proportional to their shares of the total capacity $S_j$, times a replacement frequency $T_{ij}^{-1}$ (the inverse of the changeover time constant $T_{ij}$). This is a special form of the Lotka-Volterra equation, which was used in previous work to build a model to forecast changes in the technology composition of the power sector \cite{Mercure2012}. It reduces to a logistic form in the case of two interacting technologies.

The process of changes in the number of units, following demand, efficiency or composition changes is described by processes related to the construction and decommission of units. We thus present a detailed `first principles' derivation of a theory representing this process here. The following section provides the appropriate grounding for deriving eq.~\ref{eq:LotkaVolterra}. In most previous uses for technology forecasting, eq.~\ref{eq:LotkaVolterra}. is expressed as fractions of the total capacity, and in this form, produces logistic transitions for some special circumstances. Its proportionality term, the frequency frequency $T_{ij}^{-1}$, and the absolute capacity scaling (the time dependent factor $U_{tot}$ that is divided out), are generally left ambiguous and thus open to question. This can be resolved if the equation is derived in absolute capacity terms, provided here. This exercise is helpful for better understanding the properties of this model of technology substitution, its absolute parametrisation and time scaling, and makes its quantitative application straightforward, in a way that is not empirical. In other words, time constants can be obtained without the use of fits of logistic functions to data, which in any case, as shown in this work, cover only a small subset of possible situations, which in some cases or sectors may never have actually occurred in history.\footnote{For instance, finding the diffusion time constant of solar photovoltaic devices and other new technologies from empirical fits of logistic functions is difficult since their diffusion is very preliminary, and one doesn't know $a$ $priori$ which and/or how many technologies it is replacing, and therefore, if a logistic transition is expected at all. As shown below, time constants belong to interacting pairs of technologies, not to particular technologies.}

The description of the model builds upon the property of equation~\ref{eq:LotkaVolterra} to generate an exact equivalence between an exhaustive series of of pairwise comparisons of technology options and a simultaneous comparison of all technologies, for representing investor behaviour. This property stems from its differential nature, which would not be true were the equations integrated for particular cases (as is done in \cite{Marchetti1978}, the equations of which are only true if the situation is driven by a sequence of pairwise technology transitions). A general solution does not exist; instead several possible properties emerge which include chaotic behaviour and complex dynamic outcomes.

\subsection{Definitions and assumptions}
Considering a sector where $N$ technologies compete in the marketplace and produce a single substitutable good, the system is viewed in terms of firms putting forward construction bids, which may or may not be chosen by investors or consumers. Unsuccessful bids are discarded.\footnote{This can also be viewed as items put on sale, but not necessarily bought.} I thus define the following:

\begin{description}
\item[1- The building capacity of technology $i$.] This corresponds to the number of units that could be built in a unit of time without restrictions related to how many could be sold, thus the fastest production rate possible. The assumption taken is the following: during an interval of time $t_i$, the industry is able to increase its production capacity for technology $i$ by a factor of between $e = 2.71$. However in the long run, when sales decrease the industry also scales down its activities and reduces production capacity by the same factor. Its size is therefore always proportional to $U_i$ the number of units of $i$ in operation, which is assumed to stem from a certain fraction of income reinvested. The number of units that can be built during time interval $\Delta t$ is:
\beq
{U_i \over t_i} \Delta t.
\eeq

\item[2- The total building capacity.] This is the sum over all technology categories of the previous term:
\beq
\sum_k {U_k \over t_k} \Delta t = {U_{tot}\over \overline{t}}\Delta t
\eeq
where $\overline{t}$ corresponds to the weighed average production capacity expansion time constant:
\beq
\overline{t} = \left(\sum_\ell {S_\ell \over t_\ell} \right)^{-1} , \quad S_\ell = {U_\ell \over U_{tot}}.
\eeq

\item[3- The number of units decomissioned.] During each time interval $\Delta t$, a number of units of technology $j$ come to the end of their working life which lasts a time $\tau_j$ and are retired (the lifetime can be either the physical lifetime or the lifetime of social habits, expanded below), either permanently or to be replaced. This number is a fraction of the number of units operating. Assuming that these units were constructed at rates uniform in time times the diffusion scale (e.g. not all at the same time),\footnote{Without this assumption, which is not enormously constraining, one has to track the life of every individual unit. This assumption involves constant average lifetimes, therefore no early scrapping.} this frequency is
\beq
{U_j \over \tau_j} \Delta t
\eeq

\item[4- The total number of units decommissioned.] This is the sum of the previous over all technologies:
\beq
\quad \sum_\ell {U_\ell \over \tau_\ell} \Delta t = {U_{tot}\over \overline{\tau}} \Delta t,
\eeq
where $\overline{\tau}$ is the weighed average technology lifetime:
\beq
\overline{\tau} = \left(\sum_\ell {S_\ell \over \tau_\ell} \right)^{-1}.
\eeq

\item[5- The choice between technologies.] In a pairwise comparison of technologies, the fraction of times investors or consumers will favour technology $i$ over $j$ is $F_{ij}$ (e.g. 60\% $i$ and 40\% $j$). This choice should be exclusive, i.e. $F_{ij} + F_{ji} = 1$ (and $F_{ii} = {1 \over 2}$). This matrix is time dependent and has the following property:
\beq
\sum_i \sum_j F_{ij} = {N^2 \over 2}
\eeq
Technology choices are probabilistic in nature and stem from considerations such as cost comparisons, technology performance, etc. This should be detailed specifically for individual sectors, as in \cite{Mercure2012}.

\item[6- Changes in total capacity required.] This stems from increases in demand that require new builds, or reductions in demand leading to permanent decommissions. This will be denoted respectively $\Delta U_{tot}^{\uparrow}$ and $\Delta U_{tot}^{\downarrow}$.

\item[Assumption 1: Limited number of building bids.] This rule stipulates that each firm building a particular technology type always puts forward no more bids for replacing other units and producing new units than it is able to produce in total with its full construction capacity.

\item[Assumption 2: Only successful bids are built.] The choice of the investors out of the options for the number of required units to be built restricts the number of successful construction bids. The choices are described by the matrix $F_{ij}$. Bids or proposals not chosen by investors or consumers are discarded.

\item[Assumption 3: No early scrapping.] Units operate until the end of their lifetime, where they are either replaced or decomissioned.\footnote{In the case of sudden high costs of operation, their capacity factors can however decrease to low values, meaning little use.}

\end{description}
%----------------------------------------------------------------------------------------------------------------------

\subsection{Total changes}

The evolution of a particular technology mix is assumed in this model to occur through either building units for supplying increases in demand, or to replace decommissions of old units at the end of their working lives, defined as lasting the average lifetimes $\tau_i$. The total number of units of any category built at any one time $\Delta U_i$ depends on how many units are decommissioned, how many of these are replaced and if they are, how many additional new units are  required:
\beq
\Delta U_i =  \sum_j \Delta U_{ij}\bigg|_{U_{tot}} + \Delta U_i^{\uparrow} - \Delta U_i^{\downarrow} ,
\label{eq:changes}
\eeq
where the first term denotes decommissions replaced by units of technology $i$ and decommissions of $i$ replaced by other technologies (exchanges),\footnote{The underscore $U_{tot}$ in eq.~\ref{eq:changes} indicate that the total capacity in maintained constant in the exchange term, keeping with a convention where changes in $U_{tot}$ appear in the second and third terms of this equation.} the second corresponds to  additional units due to increases in demand and the last term denotes permanent decommissions (i.e. not replaced).

\subsection{Exchange term}

\begin{figure}[t]
	\begin{center}
		\includegraphics[width=0.7\columnwidth]{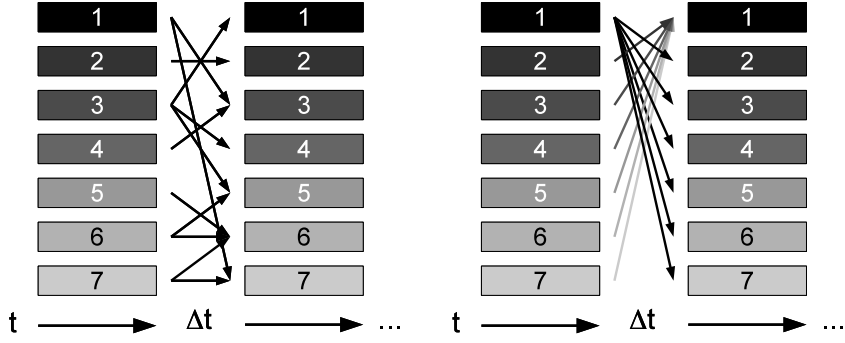}
	\end{center}
	\caption{Description of the exchange term. $Left$ Units of capacity are decommissioned and replaced, either by the same technology or by others, leading to exchanges of capacity between categories, denoted in this illustration with numbers 1-7. This gives rise to complex behaviour. $Right$ The complex exchange behaviour is simplified by exploring pairwise substitution processes individually. By summing the contribution of all losses of one technology category to all others, and its simultaneous gains at the expense of all others, one can derive the evolution of the total capacity of each technology. This illustrates eq.~\ref{eq:LVcapacity}, but can equivalently represent the shares of capacity given by eq.~\ref{eq:Shares}.}
	\label{fig:Exchanges}
\end{figure}

The first term in eq. \ref{eq:changes} corresponds to exchanges between technologies $i$ and $j$. Its structure may appear complex when first considered (see figure \ref{fig:Exchanges}, left panel), but can be decomposed in terms of individual pairwise interactions (right panel). This involves two independent terms, units decommissioned in $j$ replaced by $i$, $U_{j\rightarrow i}$ and the reverse $U_{i\rightarrow j}$ (where both can be non-zero due to the probabilistic nature of investor choices, see \cite{Mercure2012} for details).

Units of $j$ decommissioned can be expressed as a fraction (in brackets) of the total number of units decommissioned per unit time, which occurs at a frequency $\overline{\tau}^{-1}$:
\beq
{\left({U_j \over \tau_j} \over \sum_\ell {U_\ell \over \tau_\ell}\right)} {U_{tot} \over \overline{\tau}} \Delta t. 
\eeq
The fraction of these that investors would like to replace by technology $i$ is $F_{ij}$, thus giving
\beq
F_{ij} \overline{\tau} {U_j \over \tau_j} {\Delta t \over \overline{\tau}} .
\label{eq:decom}
\eeq

The number of construction bids put forward for technology $i$ corresponds to the maximum number of units of technology $i$ that could be built if the total construction capacity was used,
\beq
{\left({U_i \over t_i}  \over \sum_k {U_k \over t_k}\right)} {U_{tot}\over \overline{t}} \Delta t.
\label{eq:max_i}
\eeq
where $\overline{t}^{-1}$ must be greater than $\overline{\tau}^{-1}$, otherwise population declines cannot be fully replaced by the production capacity. This also means that the production capacity is never used fully solely for capacity replacements (and thus some bids are discarded).

Out of these $U_{tot}/\overline{t}$ construction bids put forward, only $U_{tot}/\overline{\tau}$ will be chosen, the remaining being discarded. In specific bidding cases however, it may happen that no bids of a particular kind have been put forward, while in some others several competing bids for a particular technology might have been proposed. Overall, the choice occurring in every single case must take account of the fractional composition of the total amount of bids (term in brackets in expr.~\ref{eq:max_i}) in order to determine which types of bids are available, even given the preferred choices of investors. Therefore the number of decommissions of technology $j$ that will actually be replaced by technology $i$ corresponds to expr.~\ref{eq:decom} times the fractional composition of the bidding process:
\beq
\Delta U_{j \rightarrow i} = {\left({U_i \over t_i}  \over \sum_k {U_k \over t_k}\right)} \bigg(F_{ij}\bigg) {\left({U_j \over \tau_j} \over \sum_\ell {U_\ell \over \tau_\ell}\right)}{U_{tot} \over \overline{\tau}} \Delta t = \overline{t} {U_i \over t_i} F_{ij} \overline{\tau}{U_j \over \tau_j} {1 \over U_{tot}}{\Delta t \over \overline{\tau}},
\eeq
and the remaining unsuccessful bids are discarded (analogous to newborns dying in biological systems). Therefore, the exchange term, including flows both ways, is
\beq
\Delta U_{ij} = \Delta U_{j \rightarrow i} - \Delta U_{i \rightarrow j} = \left({\overline{t} \overline{\tau} \over t_i \tau_j} F_{ij} - {\overline{t} \overline{\tau} \over t_j \tau_i} F_{ji} \right)U_i U_j {1\over U_{tot}} {\Delta t\over \overline{\tau}},
\eeq
where neither $\Delta U_{j \rightarrow i}$ nor $\Delta U_{i \rightarrow j}$ is zero, due to the probabilistic nature of investor choices (i.e. $F_{ij}$ is generally not zero, and $F_{ij} + F_{ji} = 1$). If $F_{ij} = 1/2$ then flows in both ways are equal and cancel out, but otherwise one technology category will grow at the expense of the other.

This equation couples changes between technology categories, with a time scaling constant $\overline{\tau}$, which determines at which average rate units require replacement, and a unitless substitution (or routing) matrix $A'_{ij}$, or as a matrix of substitution frequencies $A_{ij}$ (cancelling both factors $\overline{\tau}$ in the equation)
\beq
A'_{ij} = {\overline{t} \overline{\tau} \over t_i \tau_j}\quad \text{or} \quad A_{ij} = {\overline{t} \over t_i \tau_j},
\eeq
The total number of units of $i$ built during $\Delta t$ is the sum over all technologies it replaces:
\beq
\Delta U_i = \sum_j \left({\overline{t} \overline{\tau} \over t_i \tau_j} F_{ij} - {\overline{t} \overline{\tau} \over t_j \tau_i} F_{ji} \right) U_i U_j {1\over U_{tot}}{\Delta t\over \overline{\tau}}.
\label{eq:LVcapacity}
\eeq

\subsection{Increases in demand}

Changes in capacity can occur irrespective of decommissions, where bids for technologies are chosen for building a larger number of producing units in order to respond to increasing demand. The number of bids successful will depend on pairwise choices, where investors choose how to allocate $\Delta U_i^{\uparrow}$ between bids from technologies $i$ and $j$. Investor choices similarly follow the probability matrix $F_{ij}$. The number of successful bids for $i$ that were unsuccessful for $j$ is:
\beq
\Delta U_i^{\uparrow} = {1 \over N }\sum_i {\left({U_i \over \tau_i}  \over  \sum_\ell {U_\ell \over \tau_\ell}\right)}F_{ij} \Delta U_{tot}^{\uparrow} = \sum_j  {\overline{t} \over t_i}{F_{ij} \over N} U_i { \Delta U_{tot}^{\uparrow} \over U_{tot}}.
\eeq
This can be seen as first allocating equally the new builds between all technologies, $as$ $if$ they were additional decommissions to be replaced, and then following the previous reasoning for technology exchanges. If ${\Delta U_{tot}^{\uparrow} \over \Delta t} > {U_{tot} \over \overline{t}}$, then the demand exceeds the total production capacity, and cannot be met. Thus the increase in $U_{tot}$ must have a ceiling at the value of $U_{tot}/\overline{t}$ which is limited by production capacity.

\subsection{Decreases in demand}

The last term in eq.~\ref{eq:changes} corresponds to decommissions not replaced, which can also be expressed as the fraction of total units permanently decommissioned, due to demand reductions,  that belong to $i$:
\beq
\Delta U_i^{\downarrow} = {1 \over N } \sum_j F_{ji}{\left({U_i \over \tau_i} \over \sum_\ell {U_\ell \over \tau_\ell}\right)} \Delta U_{tot}^{\downarrow} = \sum_j{F_{ji}\over N}{\overline{\tau}\over \tau_i} U_i{\Delta U_{tot}^{\downarrow}\over U_{tot}},
\eeq
where units that are permanently decommissioned, out of total of $U_{tot}/\overline{\tau}$ coming to the end of their working life, are chosen as those not replaced (those being replaced belonging to the exchange term), and therefore involves the factor $F_{ji}$ (as opposed to $F_{ij}$). If however ${\Delta U_{tot}^{\downarrow} \over \Delta t} > {U_{tot}\over \overline{\tau}}$, then the decreases in demand exceed the rate at which units come to the end of their working lifetime, and this results in early scrapping (excluded here) or in reduced full load hours.

\subsection{Total changes}

The total change in technology $i$ is obtained by replacing all terms in eq. \ref{eq:changes}:
\beq
\Delta U_i = \sum_j \left[\left({\overline{t} \overline{\tau} \over t_i \tau_j} F_{ij} - {\overline{t} \overline{\tau} \over t_j \tau_i} F_{ji} \right) U_i U_j {1 \over U_{tot}}{1\over \overline{\tau}}\Delta t +  {\overline{t} \over t_i}{F_{ij} \over N} { U_i \over U_{tot}} \Delta U_{tot}^{\downarrow} - {\overline{\tau}\over \tau_i}{F_{ji}\over N} {U_i\over U_{tot}} \Delta U_{tot}^{\uparrow} \right] 
\label{eq:Total}
\eeq
The shares equation is derived from the chain derivative of $S_i(t) = U_i(t)/U_{tot}(t)$,
\beq
{\Delta S_i \over \Delta t} = {1 \over U_{tot}}\left({\Delta U_i \over \Delta t}\right) - {U_i \over U_{tot}^2}\left({\Delta U_{tot} \over \Delta t}\right).
\eeq
The total capacity can either increase or decrease, and therefore either $\Delta U_{tot}^{\uparrow}$ or $\Delta U_{tot}^{\downarrow}$ is null. In the most commonly expected case of increases in total capacity, this results in
\beq
\Delta S_i = \sum_j \left[{1 \over \overline{\tau}}\left({\overline{t} \overline{\tau} \over t_i \tau_j} F_{ij} - {\overline{t} \overline{\tau} \over t_j \tau_i} F_{ji} \right) S_i  S_j +    \left({\overline{t} \over t_i}{F_{ij} \over N} - {1\over N} \right)S_i{1 \over U_{tot}}{\Delta U_{tot}\over \Delta t}\right]\Delta t,
\label{eq:RealShares1}
\eeq
while in the case where $U_{tot}$ decreases,
\beq
\Delta S_i = \sum_j \left[{1 \over \overline{\tau}} \left({\overline{t} \overline{\tau} \over t_i \tau_j} F_{ij} - {\overline{t} \overline{\tau} \over t_j \tau_i} F_{ji} \right) S_i  S_j -  \left({\overline{\tau} \over \tau_i}{F_{ji}\over N} - {1\over N} \right)S_i{1 \over U_{tot}}{\Delta U_{tot}\over \Delta t}\right]\Delta t.
\label{eq:RealShares2}
\eeq

For specific cases where 
\beq
{\overline{t} \over t_i}\sum_j {F_{ij} \over N} = 1 \quad \text{and} \quad {\overline{\tau} \over \tau_i}\sum_j{F_{ji}\over N} = 1,
\eeq
the second term in both equations cancel exactly, leaving 
\beq
\Delta S_i = {1 \over \overline{\tau}}\sum_j \left({\overline{t} \overline{\tau} \over t_i \tau_j} F_{ij} - {\overline{t} \overline{\tau} \over t_j \tau_i} F_{ji} \right) S_i S_j \Delta t,
\label{eq:Shares}
\eeq
which is the same set of coupled differential equations of belonging to the Lotka-Volterra family as eq.~\ref{eq:LotkaVolterra}, and is the equation used in our previous work \citep{Mercure2012}. The error generated by this approximation is of the order of 
\beq
1 - {\overline{t} \over t_i}\sum_j {F_{ij} \over N}{\Delta U_{tot}\over U_{tot}} \quad \text{or} \quad 1 - {\overline{\tau} \over \tau_i}\sum_j{F_{ij}\over N}{\Delta U_{tot}\over U_{tot}}.
\eeq
This corresponds to errors in the allocation of new ($U_{tot}$ increases) or decommissioned ($U_{tot}$ decreases) units. Actually, eq. \ref{eq:Shares} would be realised if there were no restrictions on the rates of building and decommissions (i.e. if the building rate could match increasing demand exactly, and if early scrapping was allowed when capital is rapidly becoming superfluous). However, these errors in the allocation of new and scrapped units is very secondary due to the exchange term, which always enables units to change category. Therefore the difference is actually $very$ small. This is why even for large total capacity increases, the Lotka-Volterra equation expressed in fractions of the total capacity still works outstandingly well with historical data. It is  stressed that the major part of constructions at any one time aims at replacing existing units coming to the end of their working life.

In absolute capacity space, the shares equation of the FTT power sector model \citep{Mercure2012} is
\beq
\Delta U_i = \left[\sum_j \left({\overline{t} \overline{\tau} \over t_i \tau_j} F_{ij} - {\overline{t} \overline{\tau} \over t_j \tau_i} F_{ji} \right) S_i S_j {U_{tot } \over \overline{\tau}} + S_i{\Delta U_{tot} \over \Delta t}\right]\Delta t,
\eeq
which corresponds to allocating new units (i.e. those that are not replacements) or scrapped units (those not replaced), according to market shares, at equal rates between technologies, while eqs.~\ref{eq:RealShares1}-\ref{eq:RealShares2} distribute new units according to investor choices and rates of building (how much $t_i$ differs from $\overline{t}$). The difference with this equation is that technologies with faster rates of building could take up faster the new building work for increases in demand than those with slower rates, changing the allocation. The second term in eqs. \ref{eq:RealShares1}-\ref{eq:RealShares2} provide a small correction for this.

\subsection{Special case: substitutions in pairs}

A special case of the system exists, and often occurs, when only two technologies are substituted. The result, derived here, is the classic logistic transition often observed in history (see for instance \cite{Fisher1971,Bass1969,Nakicenovic1986,Sharif1976,Grubler1999,Wilson2009,Farrell1993,Marchetti1978}):
\bea
\Delta S_1 &=& \left({\overline{t} \overline{\tau} \over t_1 \tau_2} F_{12} - {\overline{t} \overline{\tau} \over t_2 \tau_1} F_{21}\right)S_1 S_2{\Delta t\over \overline{\tau}},\nn\\
&=& {1 \over T_{12}} S_1\left(1- S_1\right) \Delta t,
\eea
which, in the small time step limit, has the solution
\beq
S_1(t) = {1  \over 1 + \exp\left({-t/ T_{12}}\right)}, \quad \text{time constant} \; T_{12} = \left({\overline{t}  \over t_1 \tau_2} F_{12} - {\overline{t} \over t_2 \tau_1} F_{21}\right)^{-1}.
\eeq
This expression evolves towards the domination, in the long run, of technology 1 if the transition time constant $T_{12}$ is positive and towards 2 if it is negative. In the particular case where the choice of investors is very strongly biased, for instance towards technology 1 (2), then $F_{12} = 1$ and $F_{21} = 0$ ($F_{12} = 0$ and $F_{21} = 1$) and the time constant becomes simply $t_1 \tau_2 \over {\overline{t}}$ (${- t_2 \tau_1 \over \overline{t} }$). However, in most cases the time constant, as extracted from data, will be somewhere in the range between ${\overline{t}  \over t_1 \tau_2}$ or $-{\overline{t} \over t_2 \tau_1}$ and $\pm$ infinity.\footnote{In cases where technologies are very similar but one is slightly preferred, the transition time can be very long.} It is therefore difficult to interpret time constants obtained from fitting logistic curves to past data without additional knowledge on investor preferences $F_{ij}$.

\subsection{Special cases: limits}

\begin{figure}[t]
	\begin{center}
		\includegraphics[width=0.7\columnwidth]{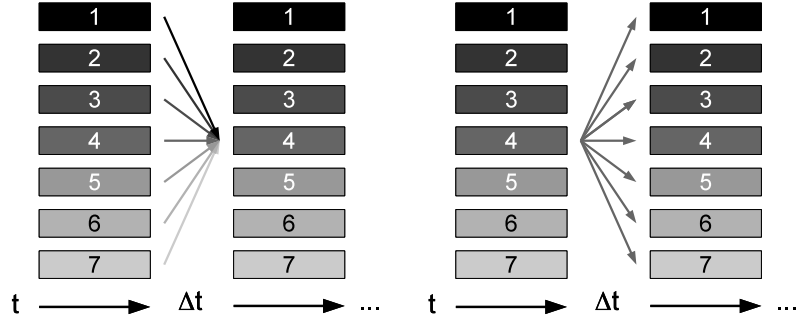}
	\end{center}
	\caption{$Left$ First limiting case where one technology is favoured over all others. $Right$ Second limiting case where one technology is disfavoured compared to all others.}
	\label{fig:Limiting}
\end{figure}

In the multi-technology interaction case where the choice between $i$ and all possible $j$ favours overwhelmingly $i$, $F_{ij} = 1$, $F_{ji} = 0$ and $F_{ii} = 1/2$ (see figure \ref{fig:Limiting}, left panel). In this case, the number of new units of $i$ is (including self replacements)
\beq
\Delta S_i = {\overline{t} \overline{\tau} \over t_i}S_i \left( \sum_j {S_j \over \tau_j} - {2 S_i F_{ii} \over \tau_i} \right) {\Delta t \over  \overline{\tau}} = {\overline{t}\over t_i} S_i \left(1- S_i {\overline{\tau}\over \tau_i}\right) {\Delta t\over \overline{\tau}},
\eeq
which in the small time step limit is an ordinary differential equation with a logistic function as a solution
\beq
S_i(t) = {1 \over  {\overline{\tau}\over \tau_i} + \exp\left({-{\overline{t}\over t_i\overline{\tau}}t}\right)}, \quad \text{time constant} \; {t_i\over \overline{t} }\overline{\tau}
\eeq
This corresponds to the highest growth rate allowed by the rate of change of the construction capacity, and thus a situation of growth limited by the construction capacity and the average decommission rate. \footnote{Note that this expression appears to converge to a value higher than 1 if $\tau_i > \overline{\tau}$, seemingly a contradiction. However, as technology $i$ comes to dominate, $\overline{\tau}$ gradually approaches $\tau_i$, and $S_i$ does not exceed 1. Therefore this example, taken for illustration, is an approximation and one must keep in mind that $\overline{\tau}$ is time dependent, even though it was assumed constant when solving the differential equation. }

In the converse case where some or all $j$ are overwhelmingly favoured over $i$,  $F_{ij} = 0$ and $F_{ji} = 1$ (see figure \ref{fig:Limiting}, right panel). The change in the number of units of $i$ becomes (including self replacements)
\beq
\Delta S_i = - {\overline{t} \overline{\tau} \over \tau_i} S_i \left( \sum_j {S_j \over t_j} - {2 S_i F_{ii} \over t_i}\right){\Delta t\over \overline{\tau}} = - {\overline{\tau}\over \tau_i}S_i \left( 1- S_i {\overline{t}\over t_i}\right) {\Delta t \over \overline{\tau}} 
\eeq
a differential equation with a logistic solution:
\beq
S_i(t) = {1  \over {\overline{t}\over t_i} + \exp\left({t/ \tau_i}\right)}, \quad \text{time constant} \; \tau_i,
\eeq
a transition that occurs at the fastest reduction rate possible for $i$ (the full decommission rate). Note that the inclusion of the self-replacement term makes the behaviour logistic, otherwise it would be an exponential decay.

These are cases limited by the fastest possible growth and decline of capacity. Therefore, all other time constants for exchanges of capacity that may be observed in any situation must lie at or above the values of $\tau_i$ (in the case of decline) and ${t_i \over \overline{t}}\overline{\tau}$ (in the cases of growth), and can involve any number of any possible combinations of terms with different $t_i$ and $\tau_i$. Actually, except for the limiting cases and single pairwise interactions, all other results will $not$ be logistic. Logistic transitions occur only in these very special circumstances.\footnote{Other cases, involving three-body or many-body interactions, cannot be solved analytically but are readily calculated numerically.} Therefore, all empirical work where a logistic transition was observed indicate special situations corresponding to either of these three possibilities (pairwise substitutions or the two limiting cases), $not$ the general case.

\subsection{Time constants and their meaning}

%Additionally, it must be noted that $\overline{t}$ and $\overline{\tau}$ are slightly time dependent. This can be calculated easily in any simulation run at every time step. However, since the lifetimes and construction rates do not differ by large amounts, they are both nearly constant. In the case of FTT:Power for the power sector, the values of $t_i$ were taken as proportional to the construction times (or lead times). In this case, $\overline{t}$ varies between 3 and 5 years, with an average near 4 years, and $\overline{\tau}$ varies between 42 and 52 years, with an average near 47 years. 

In order to scale the absolute time dependence of eq. \ref{eq:Shares}, one is required to completely resolve  any ambiguity regarding the meaning of all time constants. Time constants are taken to mean technology unit lifetimes ($\tau_i$) or the lifetime of the social habit or production method,\footnote{For some technologies, unit lifetimes are irrelevant while the time to change a method or a habit can be long, for example in communications systems, computers, software, entertainment equipment, etc, where technology changes are much slower and constrained by other considerations than unit lifetimes, or where unit lifetimes are undefined (e.g. software).}  and the time required for increasing the building capacity ($t_i$) by a factor $e$. Left as they are, these equations result in the capacity changing by at most $1\over e$ over the duration of a time constant (i.e. they are decay rates). If $\tau_i$ is assumed to correspond to a lifetime, this decay rate is too long. $\tau_i$ corresponds to the fastest possible decommission rate, in a case where no units are replaced. Thus the decay from a large population to zero should last one lifetime if units have been built up to the point where the decision was taken to stop replacing them. I thus postulate that the fastest possible rate of change for which the capacity can possibly change from 99\% to 1\% is one lifetime, corresponding to an additional factor of 10 in the time scaling.\footnote{A variation from 50\% to 99\% requires five time constants, akin for the logistic function to sigmas for gaussian functions} This is justified by the fact that, as in an extinction process occurring through unsuccessful births, it takes one lifetime to extinguish one species.

It is difficult to know exactly what the time constants are for the expansion of building capacity, as such values are not readily available, as opposed to unit lifetimes. However it can be argued that, in a similar way that learning-by-doing occurs every time the building of a unit of technology is completed, due to the rates of cash flows associated to sales, the expansion of building capacity occurs every time a new unit of technology is built and sold. Building times vary, and therefore building capacity expansion rates also vary, proportionally to the building rate. The proportionality factor is not known, however the constants $t_i$ always appear in the equations as a ratio with $\overline{t}$, and this unknown proportionality factor $always$ cancels. Therefore, for systems with long lead times (e.g. in the power sector), it is equivalent to use building time constants instead of building capacity expansion time constants. For systems with short lead times (e.g. end-use appliances), this argument is less robust and time constants for industrial production expansion should be used.

For FTT:Power, power sector capacity lifetimes vary between 20 and 80~years, with a current global average of about $\overline{\tau} = 50$~years. Building times (lead times) vary between 1 and 7~years, with an average near $\overline{t} = 4$~years.\footnote{Note that $\overline{\tau}$ and $\overline{t}$ are not exactly but only approximately constant, since the $S_i$ are time dependent. In a specific sector, values for $\tau_i$ and $t_i$ do not differ by large factors. Nevertheless, in a situation where technology moves towards systems with shorter lifetimes, the overall possible rates of technological change increase (e.g. the power sector moving towards shorter lived renewables such as wind and solar energy).} Therefore, the final time scaling of the FTT:Power shares equation is 
\beq
A_{ij} = {10 \overline{t} \over t_i \tau_j} \simeq {40 \over t_i \tau_j}.
\eeq

Changeover time constants for energy systems are known to range between 50 and 100 years. In order to to avoid any possible ambiguity regarding the resulting changeover time constants predicted by this model, we stress that the changeover time is not often near the fastest possible changeover time of one lifetime (except in the event of strong policy). Rather, it is expressed by the following term of the shares equation:
\beq
\left( F_{ij}A_{ij} - F_{ji}A_{ji}\right) \simeq \left(F_{ij}{40\over t_i \tau_j} - F_{ji} {40 \over t_j \tau_i}\right).
\eeq

\subsection{Constraints of the market that reflect imperfect competition}

\begin{figure}[t]
	\begin{center}
		\includegraphics[width=0.7\columnwidth]{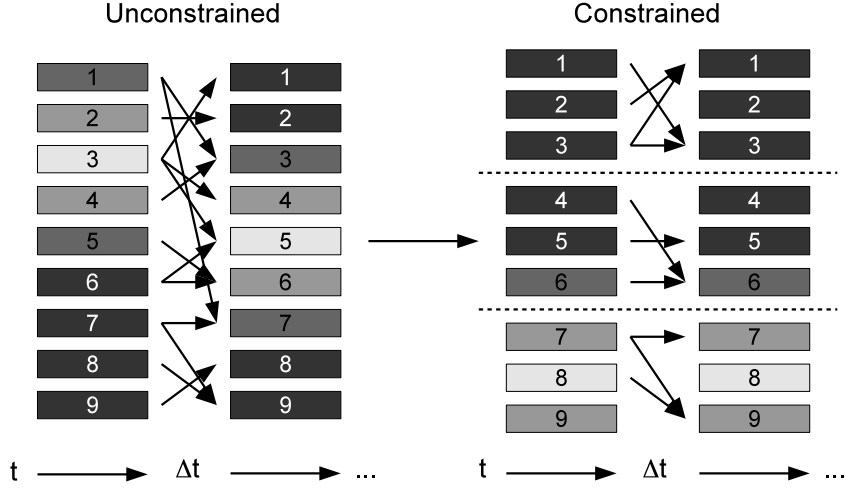}
	\end{center}
	\caption{Market splitting that arises with market constraints. $Left$ In an unconstrained market, substitution processes can occur between all possible pairs of systems, competing on a single price level. $Right$ In a constrained system, the market segments and substitutions occur between technologies within submarkets with different price levels. The cumulative share within each submarket is maintained, and three constraints are used for illustration purpose.}
	\label{fig:MarketSplitting}
\end{figure}

Technologies were assumed in the preceding as perfectly substitutable for one another, which implies that the service that they produce, within a category, is identical and independent of the particular technology used. This is of course too big a simplification of reality, which could seemingly invalidate the model. The market may have requirements that constrains the shares of particular types of systems. Such constraints must restrict the number of possible technology mixes that are possible. 

This was discussed at length for the power sector in previous work in terms of three types of electricity generation technologies and two types of constraints \citep{Mercure2012}. The magnitude of peak load demand, as well as the amount of variable renewable capacity, requires a certain amount of flexible electricity generation systems, unless significant amounts of energy storage is available. This reflects an engineering grid stability constraint. I had argued there that coupled constraints can be introduced, related to investor reticence of investing in assets that could become unused or stranded, and therefore investment in particular types of technologies slows down when the limit of stability of the grid is approached. In situations limited by flexibility, limits applied to market shares generate a spontaneous splitting of the market, where peak load technologies, variable systems and base load generation systems interact in respective separate submarkets at different price levels.

This can be generalised for any type of sets of competing technologies, for example for transport, lighting, etc. The demand generally has specific inflexible characteristics which can segregate the market into segments. For example, transport demand can be classified in terms of distance travelled, and certain types of substitution between technologies serving different demand segments are unlikely (e.g. motorcycles for long distance travelling, or using airplanes for commuting, etc.). There is, however, a certain amount of overlap between demand segments, which does not allows them to be treated independently (e.g. substitution can occur between air and train travel, as well as between car and train travel, but to a lesser extent between car and air travel.) This reflects the variety of motivations that generate the aggregate transport demand. Urban design and the availability of parking space influence transport choices.

Details underlying the structure of the demand for a service cannot be completely detailed in a model of technological change, but must be summarised appropriately in order to avoid generating unlikely or unfeasible technology mixes. For this, dynamic limits to market shares are introduced. We define $n$ segments of service demand or constraints, expressed either in units of service generation, i.e. $D_1$, $D_2$, ... $D_n$,  or in units of capacity $U_1$, $U_2$, ... $U_n$. As many share limit $\hat{S}_i$ equations can be defined such as 
\bea
\hat{S}_i^1 = S_i \pm \left[ \sum_j a_j S_j - {\overline{CF}\over CF_i}{D_1\over D}\right], && 
\hat{S}_i^1 = S_i \pm \left[ \sum_j a_j S_j - {U_1\over U_{tot}}\right],\nn\\
\hat{S}_i^2 = S_i \pm \left[ \sum_j a_j S_j - {\overline{CF}\over CF_i}{D_2\over D}\right], &\text{or}&
\hat{S}_i^2 = S_i \pm \left[ \sum_j a_j S_j - {U_2\over U_{tot}}\right],\nn\\
..., \quad\quad&&\quad\quad ...,\nn\\
\hat{S}_i^n = S_i \pm \left[ \sum_j a_j S_j - {\overline{CF}\over CF_i}{D_n\over D}\right], &&
\hat{S}_i^n = S_i \pm \left[ \sum_j a_j S_j - {U_n\over U_{tot}}\right], 
\eea
where $\hat{S}_i$ denotes a limit that the share value of technology $i$, $S_i$, is not allowed to cross, and $a_i$ sets the type of service contribution that technology $i$ provides, and is either equal to 1 or -1. This can be an upper limit (positive sign) or a lower limit (negative sign). The most constraining limit is then taken. These limits depict the $maximum$ $additional$ contribution that a technology could provide given the $current$ contributions produced by every other technology.

In such a system of share limits, when all the share values $S_i$ are far from their respective limits, exchanges can occur between every type of technology and every other type. However, when some of the $S_i$ belonging to a group of contribution that responds to a particular type of demand $D_j$ come close to their limit, they may still undergo substitutions within their group but not anymore with technologies outside this group (for example substitution between peak load technologies in the power sector in a situation limited in flexibility, or substitution between long distance transport modes). In this case a submarket associated with this demand segment $D_j$ emerges that can operate on a different cost level than the rest of the market. This is explained by invoking the argument that technology substitution across segments is suppressed by investor reticence in seeing their assets stranded, thus restricting their choices further than dictated by the matrix $F_{ij}$. This is depicted in figure \ref{fig:MarketSplitting}, where an unconstrained substitution situation is depicted on the left, and a spontaneous market splitting is described on the right, using three market segments. Such dynamic constraints can produce synergies and interference between technologies, where for instance, wind turbines can proliferate in tandem with gas turbines, but a massive expansion in gas turbines due to a subsidy for wind power can block the market for other renewables such as solid biomass based electricity.

%------------------------------------------------------------------------------------------------------------------------
%------------------------------------------------------------------------------------------------------------------------

\section{Discussion: price or policy induced efficiency changes}
\subsection{Revisiting the technology ladder}

\begin{figure}[t]
	\begin{center}
		\includegraphics[width=0.7\columnwidth]{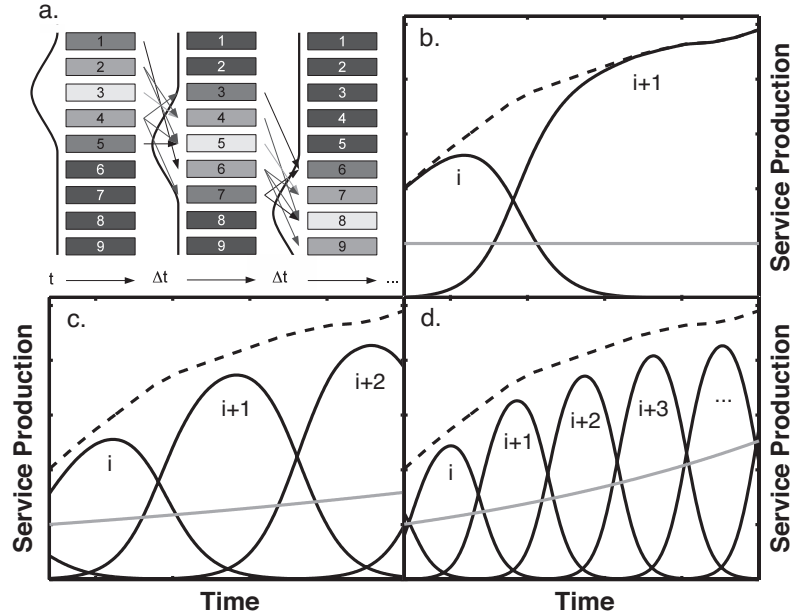}
	\end{center}
	\caption{$a.$ Process of the gradual substitution of technology towards ever more efficient systems. $b.$ Depiction of a logistic transition, where the demand for a service is given with a dashed line and the contribution of each technology is given with solid lines. $c.$ Depiction of the principle of the technology ladder, where technology $i$ is replaced by technology $i+1$, which is in turn replaced by technology $i+2$ etc. $d.$ Faster rate of change in the technology ladder. }
	\label{fig:LadderTechFor2}
\end{figure}

The concept of the technology ladder was introduced in previous work \cite{Mercure2012}. It depicts a process of gradual technology exchanges towards ever more efficient technologies given a technological change driver based onto penalising inefficiency (e.g. fuel taxes, efficiency standards, the carbon price, etc, for a history of effective technology standards in the Japanese automobile industry, see \cite{Sano2011}).\footnote{In the case of power technologies, the efficiency was defined in \cite{Mercure2012} in terms of CO$_2$ per unit of electricity produced. In the general case, several definitions could be given, e.g. gCO$_2$/km, GJ/km, GJ/kg, etc.} This is depicted schematically in panel $a.$ of figure~\ref{fig:LadderTechFor2}. For a particular type of technology providing a specific service (e.g. transport), a technology list can be drawn, with associated population, or market shares, of each component. The process of research and development continuously adds new elements to this list, while some older models are phased out. 

The drivers of change, related but not exclusively associated to technology costs, are dynamic. In an environment favouring improved efficiency, investors or consumers tend to go down the list towards more efficient technologies following the incentives given (reflected in $F_{ij}$). Due to learning by doing, however, new models tend to incur higher costs to consumers than old ones, limiting this process, but these costs decrease with cumulative production following learning curves. Therefore, the newest components of the list are difficult to access, and their diffusion beyond niche markets is obstructed by high prices. Consumers thus tend to choose on average accessible intermediate solutions, but small niches exist, as well as early adopters who might find high prices acceptable. Niches and early adopters enable to build capacity despite high prices, bringing new technology prices down along their learning curves. Thus costs of new technologies gradually decrease, and consumers gradually make their way towards ever more efficient technologies, as they buy and later discard units. The panel $a.$ of figure~\ref{fig:LadderTechFor2} depicts this for different times with a numbered list of technologies, in order of increasing efficiency. The population of the list is depicted on the left with a curve, and technology substitutions are indicated with arrows.  

The components of the list of technologies are substituted through the process described in the first part of this work, and therefore follow eqs.~\ref{eq:RealShares1}-\ref{eq:RealShares2}. In the technology ladder described here, exchanges can be assumed to follow approximately a series of logistic transitions between every consecutive pair of members of the list. For a constant value for the driver, the situation observed is likely to be that depicted in panel~$b.$ of figure~\ref{fig:LadderTechFor2},  for which a technology transition may slowly occur, where the demand for the service is given with a dashed line and the contribution of each technology is given with solid lines. However, with a dynamic driver based onto a gradually changing requirement for improved efficiency (e.g. as occurred in Japan, for instance, in the automobile industry \cite{Sano2011}), the situation depicted in panel~$c.$ may be observed, where members of the list of technologies gradually come in vogue before becoming phased out again, replaced by another ever more efficient technology. This process may be made to accelerate with a stronger rate of change in the driver, shown in panel~$d.$. However, this cannot be made to progress faster than the fastest rate of decommission allowed by the lifetime values of the technologies, unless they are scrapped early. In this limiting situation, every technology appears briefly in the marketplace for the length of one lifetime.

This process is irreversible. This is due to the learning process, which, assuming no loss of knowledge,\footnote{which can actually happen if a type of system becomes disused for too long, which in this work is considered indistinguishable from a share value of nearly zero (e.g. knowledge loss for steam locomotives). Systems phased out cannot reappear.} is itself irreversible, changing the technology choice landscape permanently. Technology choices do not tend to go backwards but only forwards in the direction of technology improvements.  The rate of change, however, is not predetermined but depends onto market conditions, making this phenomenon highly path dependent. In other words, technological change can happen in a myriad of ways. As discussed earlier, such path dependence is associated with hysteresis and irreversibility. For example, a strong sudden change in environment (e.g. an oil shock) can generate faster changes in technology properties (e.g. internal combustion engine efficiency), which (1) do not revert back once the environment returns to it's original state but remain and (2) leave permanent cost reductions produced by learning. The property of path dependence, hysteresis and irreversibility prohibits the use of a historical trend for forecasting future changes in technology properties, most notably the efficiency of energy use (such as in the AEEI parameter), since the efficiency can not be expressed as a single valued function of other variables. In particular, the rate of efficiency change generated by strong climate policy cannot be reliably predicted using historical data and without involving trends related to technology diffusion processes.

\subsection{Evaluating policy or price induced efficiency changes}

As seen with eq.~\ref{eq:Shares}, substitutions between technologies of a list may be expressed in terms of shares of production capacity. For a list of technologies with identifier $i$, ordered with increasing efficiency, any average quantity of the system can be calculated using a sum weighted by shares. For example, the average efficiency may be calculated from individual efficiency values $\alpha_i$. This average efficiency $\overline{\alpha}$ may be weighed by the relative population of each member of the list $S_i$:
\beq
\overline{\alpha}(t) = \sum_i S_i(t) \alpha_i.
\eeq
Since the population of each member of the list changes gradually following perhaps logistic transitions or series of these,\footnote{Changes may or may not be logistic; however equation~\ref{eq:Shares} is always well behaved and produces continuous and smooth changes.} the average efficiency follows gradual downwards changes. 

\begin{figure}[t]
	\begin{center}
		\includegraphics[angle = -90, width=1\columnwidth]{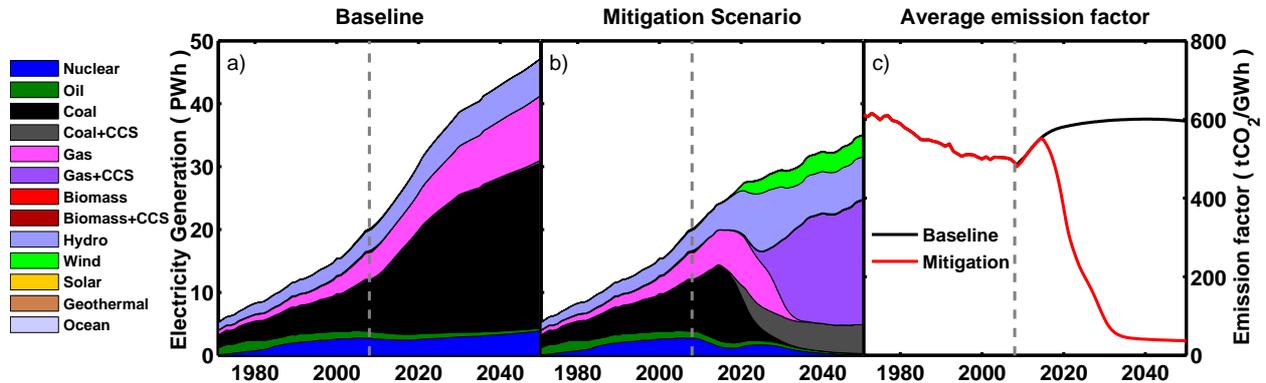}
	\end{center}
	\caption{Example of global technology forecasting scenarios for the power sector using the model FTT:Power. Historical data is shown to the left of the dashed line, while the projections are given on the right. $a)$ Baseline scenario based on current energy policy. $b)$ Climate change mitigation policy scenario (see text). $c)$ Evolution of the global average emission factor for both scenarios.}
	\label{fig:Scenarios}
\end{figure}

As an example, the power sector is taken using the model FTT:Power, shown in figure \ref{fig:Scenarios}. Two technology forecasting scenarios are given in panels $a)$ and $b)$, in terms of electricity generation, with their associated emissions factors (in tCO$_2$/GWh) given in panel $c)$. The baseline scenario assumes that current policies are maintained (based on the new policies scenario of the IEA \cite{IEAWEO2010,Mercure2012}), while the mitigation scenario assumes reduced energy demand as well as support for renewables (partly based onto the 450 scenario of the IEA \cite{IEAWEO2010}), as well as a carbon price that initially increases to near 200\$/tCO$_2$ up to around 2025, and saturates there. In particular, a subsidy is given to wind power and carbon capture and storage (CCS) technologies. Data in all panels to the left of the dashed lines are historical. The emissions factor has historically gradually decreased monotonically while electricity generation has increased approximately exponentially. In the projection, however, it initially increases before stabilising (baseline) or decreasing sharply (mitigation). The increase is due to increasing energy demand and increasing use of coal, replacing, among other things, a decreasing fleet of nuclear power stations. The stabilisation in the baseline is related to the increasing use of gas instead of coal. Meanwhile, the sharp decline in the mitigation scenario is mainly related to the gradual diffusion of CCS both for coal and gas power systems. These diffuse approximately logistically. This example demonstrates how hopeless it would be to attempt to reproduce these trends in emissions factors using historical extrapolations, simply due to the fact that too many correlated processes take place simultaneously, which could not be simply predicted from past history.

This theory has important implications for attempts at projecting scenarios for the future carbon price. A wide disagreement exists in the literature as to even the order of magnitude of the carbon price necessary to enable the technological changes required to decarbonise the global economy and prevent excessive climate change beyond 2$^{\circ}$C. Prices as high as 1700\$/tCO$_2$ and as low as 100\$/tCO$_2$ have been claimed for the year 2100 (see for instance \cite{Edenhofer2006}). Values projected depend strongly on the amount of substitution possibilities included in the models used and whether learning curves are included. The price of carbon is in essence a measure of the financial and technical difficulty associated with technological change for emissions reductions. The more emissions reducing technological change options exist, the lower the price of carbon is likely to be \cite{Edenhofer2006}. These are similar to degrees of freedom. For instance, models using an AEEI have grounded part of their process of technological change independently of the price of carbon or other fuels \cite{Azar1999}. In particular, the rate of diffusion and learning-by-doing does not respond to the price of carbon, and therefore remains rigid. It is no surprise therefore that such models predict high prices of carbon. In contrast, according to FTT:Power, the global power sector can be decarbonised with a much lower price of carbon, less than 200\$/tCO$_2$ (see \cite{Mercure2012} and figure \ref{fig:Scenarios} above). It is clear that more accurate projections of the carbon price require more elaborate dynamics of technological change to be included into models, as for instance given here.

\subsection{Parallel with evolutionary dynamics}

Without changes in the environment, natural selection and evolution would not occur. Species would adapt perfectly to their unchanging environment. This is known to have occurred for some species, for instance some crustaceans (e.g. the horseshoe crab), which are considered some of the oldest non-extinct species on earth. Changes in environment generates extinction processes, where a few individuals with specific superior genetic characteristics happen to cope with the changes better than others, and may reproduce and grow in numbers and eventually entirely replace the original population. This is the process of natural selection.

In reality, the genetic bath is diverse and features myriads of small differences between individuals. Natural selection reduces this diversity through the deaths of individuals with inferior genetic characteristics, but genetic diversity is gradually replenished through genetic mutations. This diversity combined with changes in environment is key to the process of evolution. 

Evolutionary dynamics can be applied to technological change, a subject studied by Dercole and Dercole $et$ $al.$ \cite{DercoleThesis, Dercole2008}. Technology diversity is key to technological change as it provides options that enable technology transitions. For a change in environment, that may be related to changes in prices (e.g. the price of oil or electricity) or policy (e.g. efficiency standards, taxes, price of carbon), such diversity, composed of technologies that exist in small niche markets, provide solutions that may diffuse massively and replace commonly used systems. 

For example, solar panels were initially developed for space applications (i.e. powering satellites), and have been and are currently still significantly more costly than common coal power stations, when measured using for instance their levelised cost. For a certain level of carbon emissions taxing and learning-by-doing, however, they have the potential to replace to some extent current energy systems. Thus, changes in price or policy environments may uncover many niche technologies and transform existing systems, a process which is the key to, for instance, climate change mitigation. This description is consistent with current views on technology transitions in the energy sector (see for instance in Wilson and Gr\"ubler \cite{Wilson2011}, and Gr\"ubler \cite{Grubler2012}).

\section{Conclusion}

This paper presents a theory to predict timescales of technology transitions, a subject that has never been clarified before in the large empirical literature on logistic technology substitution trends. While timescales are often seen as a constant of the system belonging to specific technologies, they are presented here as related to a particular situation occurring as a special case of a multi-technology interaction. Interaction timescales are given in terms of lifetimes and lead times and investor or consumer preferences. Technological change is represented using a family of coupled non-linear ordinary differential equations, a subset of the Lotka-Volterra family of equations originally derived for projecting the evolution of species in competing biological environments. While most of the empirical literature on technology transitions focuses on fitting logistic curves to data, it is demonstrated here that logistic transitions are only special cases of a more general system of multi-technology interactions, and that only the differential form of these equations should be used in order to appropriately represent real simultaneous multi-technology interactions.

This theory is used to describe the mechanics underlying the complex phenomenon of price or policy induced efficiency changes, or changes in any other technology properties. This phenomenon is highly path dependent, displaying hysteresis and irreversibility. This property stems from the process of learning-by-doing, which is itself irreversible, giving rise to continuous changes in the landscape of technology possibilities offered to investors or consumers. This is used to depict the process of gradual improvements in efficiency given ever changing policy or price drivers. In this picture most investors or consumers tend to choose accessible intermediate solutions, while market niches and early adopters enable learning-by-doing cost reductions of new technologies to occur, before their wide diffusion can happen. This irreversible process is consistent with for instance transport efficiency changes occurring with fuel price hikes. Finally, this framework is compared with that of evolutionary dynamics applied to technological change. 

\section*{Acknowledgements}

The author would like to thank particularly T. S. Barker (4CMR) for his guidance and support, as well as P. Salas (4CMR) for continuous discussion and common modelling work, and H. Pollitt, P. Summerton for support at Cambridge Econometrics. I would also like to thank A. Gr\"ubler (IIASA) for a lively discussion from which emerged a strong incentive to write down this theory, and M. Obersteiner (IIASA) for invoking the useful parallel with evolutionary dynamics. This work was supported by the Three Guineas Trust.

\section*{References}



%% References


\bibliographystyle{elsarticle-num}
\bibliography{CamRefs}

\begin{thebibliography}{10}
\expandafter\ifx\csname url\endcsname\relax
  \def\url#1{\texttt{#1}}\fi
\expandafter\ifx\csname urlprefix\endcsname\relax\def\urlprefix{URL }\fi
\expandafter\ifx\csname href\endcsname\relax
  \def\href#1#2{#2} \def\path#1{#1}\fi

\bibitem{Wilson2011}
C.~Wilson, A.~Gr\"ubler, {Lessons from the history of technological change for
  clean energy scenarios and policies}, {Natural Resources Forum} {35}~({3,
  SI}) ({2011}) {165--184}.

\bibitem{Grubler2012}
A.~Gr\"ubler, Energy transitions research: Insights and cautionary tales,
  Energy Policy~(0) (2012) --.

\bibitem{Grubler1999}
A.~Gr\"ubler, N.~Nakicenovic, D.~Victor, Dynamics of energy technologies and
  global change, Energy Policy 27~(5) (1999) 247--280.

\bibitem{Grubler1998}
A.~Gr\"ubler, Technology and Global Change, Cambridge University Press, 1998.

\bibitem{Azar1999}
C.~Azar, H.~Dowlatabadi, {A review of technical change in assessment of climate
  policy}, {Annu. Rev. Energy Environ.} {24} ({1999}) {513--544}.

\bibitem{IPCCSRES}
IPCC, Emission Scenarios, Cambridge University Press, 2000.

\bibitem{Jones2004}
D.~W. Jones, P.~N. Leiby, I.~K. Paik, {Oil price shocks and the macroeconomy:
  What has been learned since 1996}, {Energy Journal} {25}~({2}) ({2004})
  {1--32}.

\bibitem{MercureSalas2012}
J.-F. Mercure, P.~Salas, An assessement of global energy resource economic
  potentials, Energy~(-) (2012) --.
\newblock \href {http://dx.doi.org/10.1016/j.energy.2012.08.018}
  {\path{doi:10.1016/j.energy.2012.08.018}}.

\bibitem{MercureSalas2012b}
J.-F. Mercure, P.~Salas, On the global economic potentials and marginal costs
  of non-renewable resources, In preparation.

\bibitem{Kohler2006}
J.~Koehler, M.~Grubb, D.~Popp, O.~Edenhofer, The transition to endogenous
  technical change in climate-economy models: A technical overview to the
  innovation modeling comparison project, Energy Journal~(Sp. Iss. 1) (2006)
  {17--55}.

\bibitem{Grubb2006}
M.~Grubb, C.~Carraro, J.~Schellnhuber, {Technological change for atmospheric
  stabilization: Introductory overview to the Innovation Modeling Comparison
  Project}, {Energy Journal}~(Sp. Iss. 1) ({2006}) {1--16}.

\bibitem{Crassous2006}
R.~Crassous, J.-C. Hourcade, O.~Sassi, {Endogenous structural change and
  climate targets modeling experiments with Imaclim-R}, {Energy Journal}~(Sp.
  Iss. 1) ({2006}) {259--276}.

\bibitem{Edenhofer2006}
O.~Edenhofer, K.~Lessmann, C.~Kemfert, M.~Grubb, J.~Koehler, {Induced
  technological change: Exploring its implications for the economics of
  atmospheric stabilization: Synthesis report from the Innovation Modeling
  Comparison Project}, {Energy Journal}~(Sp. Iss. 1) ({2006}) {57--107}.

\bibitem{IMCP2006}
{Innovation Modelling Comparison Project}, Endogenous Technological Change and
  the Economics of Atmospheric Stabilisation, Special Edition, Energy Journal,
  International Association for Energy Economics, 2006.

\bibitem{Kaufmann2004}
R.~Kaufmann, {The mechanisms for autonomous energy efficiency increases: A
  cointegration analysis of the US energy/GDP ratio}, {Energy Journal}
  {25}~({1}) ({2004}) {63--86}.

\bibitem{Dowlatabadi2006}
H.~Dowlatabadi, M.~A. Oravetz, {US long-term energy intensity: Backcast and
  projection}, {Energy Policy} {34}~({17}) ({2006}) {3245--3256}.

\bibitem{IPCCAR4Ch11}
T.~Barker, I.~Bashmakov, A.~Alharthi, M.~Amann, L.~Cifuentes, J.~Drexhage,
  M.~Duan, O.~Edenhofer, B.~Flannery, M.~Grubb, M.~Hoogwijk, F.~I. Ibitoye,
  C.~J. Jepma, W.~Pizer, K.~Yamaji, Mitigation from a cross-sectoral
  perspective, in: Climate Change 2007: Mitigation. Contribution of Working
  Group III to the Fourth Assessment Report of the Intergovernmental Panel on
  Climate Change, Cambridge University Press, 2007.

\bibitem{Edenhofer2010}
O.~Edenhofer, B.~Knopf, T.~Barker, L.~Baumstark, E.~Bellevrat, B.~Chateau,
  P.~Criqui, M.~Isaac, A.~Kitous, C.~Kypreos, M.~Leimbach, K.~Lessmann,
  B.~Magne, S.~Scrieciu, H.~Turton, D.~P. van Vuuren, {The Economics of Low
  Stabilization: Model Comparison of Mitigation Strategies and Costs}, {ENERGY
  JOURNAL} {31}~({1}) ({2010}) {11--48}.

\bibitem{Marchetti1978}
C.~Marchetti, N.~Nakicenovic,
  \href{http://www.iiasa.ac.at/Research/TNT/WEB/PUB/RR/rr-79-13.pdf}{The
  dynamics of energy systems and the logistic substitution model}, Tech. rep.,
  IIASA (1978).
\newline\urlprefix\url{http://www.iiasa.ac.at/Research/TNT/WEB/PUB/RR/rr-79-13.pdf}

\bibitem{Bass1969}
F.~M. Bass, {New Product Growth for Model Consumer Durables}, {Management
  Science Series A-theory} {15}~({5}) ({1969}) {215--227}.

\bibitem{Fisher1971}
J.~Fisher, R.~Pry, A simple substitution model of technological change,
  Technological Forecasting and Social Change 3~(0) (1971) 75 -- 88.

\bibitem{Nakicenovic1986}
N.~Nakicenovic, {The automobile road to technological-change - Diffusion of the
  automobile as a process of technological substitution}, {Technological
  Forecasting and Social Change} {29}~({4}) ({1986}) {309--340}.

\bibitem{Sharif1976}
M.~N. Sharif, C.~Kabir, Generalized model for forecasting technological
  substitution, Technological Forecasting and Social Change 8~(4) (1976)
  353--364.

\bibitem{Wilson2009}
C.~Wilson, Meta-analysis of unit and industry level scaling dynamics in energy
  technologies and climate change mitigation scenarios, Tech. Rep. IR-09-029,
  IIASA (2009).

\bibitem{Farrell1993}
C.~J. Farrell, A theory of technological progress, Technological Forecasting
  and Social Change 44~(2) (1993) 161 -- 178.

\bibitem{Bhargava1989}
S.~C. Bhargava, Generalized lotka-volterra equations and the mechanism of
  technological substitution, Technological Forecasting and Social Change
  35~(4) (1989) 319--326.

\bibitem{Lotka1925}
A.~J. Lotka, Elements of Physical Biology, Wiliams and Wilkins Company, 1925.

\bibitem{Volterra1939}
V.~Volterra, The general equations of biological strife in the case of
  historical actions, Proceedings of the Edinburgh Mathematical Society 6~(1)
  (1939) 4.

\bibitem{Verhulst1838}
P.-F. Verhulst, Notice sur la loi que la population poursuit dans son
  accroissement, Correspondance mathematique et physique 10 (1838) 113--121.

\bibitem{Grubler1990}
A.~Gr\"ubler, The rise and fall of infrastructures: Dynamics of evolution and
  technological change in transport, Tech. rep., IIASA (1990).

\bibitem{Grubler2012Comm}
A.~Gr\"ubler, {Private communication} (2012).

\bibitem{Mercure2012}
J.~F. Mercure, Ftt:power : A global model of the power sector with induced
  technological change and natural resource depletion, Energy Policy 48~(0)
  (2012) 799 -- 811.
\newblock \href {http://dx.doi.org/10.1016/j.enpol.2012.06.025}
  {\path{doi:10.1016/j.enpol.2012.06.025}}.

\bibitem{Cologni2009}
A.~Cologni, M.~Manera, {The asymmetric effects of oil shocks on output growth:
  A Markov-Switching analysis for the G-7 countries}, {Economic Modelling}
  {26}~({1}) ({2009}) {1--29}.

\bibitem{Dercole2008}
F.~Dercole, U.~Dieckmann, M.~Obersteiner, S.~Rinaldi, {Adaptive dynamics and
  technological change}, {Technovation} {28}~({6}) ({2008}) {335--348}.

\bibitem{Sano2011}
M.~Sano, M.~Kii, H.~Miyoshi, Automotive technology and public policy in japan:
  a historical survey, in: H.~Miyoshi, M.~Kii (Eds.), Technological innovation
  and public policy, Pargrave Macmillan, 2011, pp. 15--43.

\bibitem{IEAWEO2010}
IEA, World Energy Outlook 2010, IEA/OECD, 2010.

\bibitem{DercoleThesis}
F.~Dercole, Evolutionary dynamics through bifurcation analysis: Methods and
  applications, Ph.D. thesis, Politecnico di Milano (2002).

\end{thebibliography}






% \end{thebibliography}

\end{document}